\documentclass[12pt]{article}
\usepackage{amssymb}
\usepackage{mathrsfs}
\usepackage{pifont}
\usepackage{tipa}
\usepackage{amsmath}
\usepackage{amssymb,amsmath}
\usepackage{amsfonts}
\begin{document}
\title{\textbf{Automorphisms and Verma modules for Generalized Schr\"{o}dinger-Virasoro algebras}
\thanks{Supported by the National Natural Science Foundation of
China (No. 10671160).}}
\author{Shaobin Tan $^1$, Xiufu Zhang $^{1,2}$\\
{\scriptsize 1. School of Mathematical Sciences, Xiamen University,
Xiamen 361005,  China}\\
{\scriptsize 2. School of Mathematical Sciences, Xuzhou Normal
University, Xuzhou 221116, China}}
\date{}
\maketitle

\begin{abstract}
Let $\mathbb{F}$ be a field of characteristic 0, $G$ an additive
subgroup of $\mathbb{F}$, $\alpha\in \mathbb{F}$ satisfying
$\alpha\notin G, 2\alpha\in G$. We define a class of
infinite-dimensional Lie algebras which are called generalized
Schr\"{o}dinger-Virasoro algebras and use $\mathfrak{gsv}[G,\alpha]$
to denote the one corresponding to $G$ and $\alpha$. In this paper
the automorphism group and irreducibility of Verma modules for
$\mathfrak{gsv}[G,\alpha]$ are completely determined.
\vspace{2mm} \vspace{2mm}
\end{abstract}

\vskip 5mm

\noindent{\bf{1. Introduction}}\vspace{3mm}

\vskip 3mm The Schr\"{o}dinger-Virasoro algebra $\mathfrak{sv}$ is
defined to be a Lie algebra with $\mathbb{F}$-basis
$\{L_n,M_n,Y_{n+\frac{1}{2}}\mid n\in \mathbb{Z}\}$ subject to the
following Lie brackets:
$$[L_m,L_n]=(n-m)L_{n+m},
[L_m,M_n]=nM_{n+m},[L_m,Y_{n+\frac{1}{2}}]=(n+\frac{1-m}{2})Y_{m+n+\frac{1}{2}},$$
$$[M_m,M_n]=0, [Y_{m+\frac{1}{2}},Y_{n+\frac{1}{2}}]=(n-m)M_{m+n+1},
[M_m,Y_{n+\frac{1}{2}}]=0.$$ It is easy to see that $\mathfrak{sv}$
is a semi-direct product of the Witt algebra
$\mathfrak{Vir_0}=\bigoplus\mathbb{C}L_n$ and the two-step nilpotent
infinite-dimensional Lie algebra
$\mathfrak{h}=\bigoplus\mathbb{C}M_n\bigoplus
\mathbb{C}Y_{n+\frac{1}{2}}$. This infinite-dimensional Lie algebra
was originally introduced in [3] by looking at the invariance of the
free Schr\"{o}dinger equation in (1+1)dimensions $(2\mathcal
{M}\partial_t-\partial_r^{2})\psi=0.$ The structure and
representation theory of $\mathfrak{sv}$ have been studied by C.
Roger and J. Unterberger in [11]. The irreducible weight modules
with finite-dimensional weight spaces over $\mathfrak{sv}$ are
classified in [7].

In order to investigate vertex representations of $\mathfrak{sv}$,
J. Unterberger (see [14]) introduced a class of infinite-dimensional
Lie algebras $\mathfrak{\widetilde{sv}}$ called the extended
Schr\"{o}dinger-Virasoro Lie algebra, which can be viewed as an
extension of $\mathfrak{sv}$ by a conformal current with conformal
weight 1. The extended Schr\"{o}dinger-Virasoro Lie algebra
$\mathfrak{\widetilde{sv}}$ is a vecter space spanned by a basis
$\{L_n,M_n,N_n,Y_{n+\frac{1}{2}}\mid n\in \mathbb{Z}\}$ with the
following commutation relations: $$[L_m,L_n]=(n-m)L_{n+m},
[M_m,M_n]=0, [N_m,N_n]=0,[N_m,M_n]=2M_{m+n},$$
$$[Y_{m+\frac{1}{2}},Y_{n+\frac{1}{2}}]=(n-m)M_{m+n+1},
[L_m,M_n]=nM_{n+m},[L_m,N_n]=nN_{n+m},$$
$$[L_m,Y_{n+\frac{1}{2}}]=(n+\frac{1-m}{2})Y_{m+n+\frac{1}{2}}
[N_m,Y_{n+\frac{1}{2}}]=Y_{m+n+\frac{1}{2}},
[M_m,Y_{n+\frac{1}{2}}]=0.$$ For all $m,n\in \mathbb{Z}$. The
structure of $\widetilde{\mathfrak{sv}}$ has been studied in [2].

Recently, a number of new classes of infinite-dimensional Lie
algebras over a field of characteristic 0 were discovered by several
authors (see [6] [8] [12] [13]). Among those algebras , are the
generalized Witt algebras, the generalized Virasoro algebras, the
Lie algebras of generalized Weyl type and the generalized
Heisenberg-Virasoro algebra.

Let $\mathbb{F}$ be a field of characteristic 0, $G$ an additive
subgroup of $\mathbb{F}$, $\alpha\in \mathbb{F}$ such that
$\alpha\notin G, 2\alpha\in G$. Motivated by the above algebras, we
introduce a new class of Lie algebras which are called the
generalized Schr\"{o}dinger-Virasoro Lie algebras which include the
Schr\"{o}dinger-Virasoro Lie algebra $\mathfrak{sv}$ as a special
case. We use $\mathfrak{gsv}[G,\alpha]$ to denote the one
corresponding to $G$ and $\alpha$. In this paper we mainly study the
automorphism group $Aut(\mathfrak{gsv}[G,\alpha])$ and the
irreducibility of Verma modules over the generalized
Schr\"{o}dinger-Virasoro Lie algebra $\mathfrak{gsv}[G,\alpha]$.

The paper is organized as follows. In section 2, we introduce the
generalized Schr\"{o}dinger-Virasoro Lie algebra
$\mathfrak{gsv}[G,\alpha]$. The necessary and sufficient conditions
of isomorphism between two of these algebras are determined. In
section 3, we determine the automorphism group of
$\mathfrak{gsv}[G,\alpha]$. In section 4, a Verma module $V(c,h)$
over the generalized Schr\"{o}dinger-Virasoro Lie algebra
$\mathfrak{gsv}[G,\alpha]$ is defined and its irreducibility is
completely determined.

\vskip 3mm

Throughout the article, we denote $\mathbb{Z}$ the set of integers,
$\mathbb{N}$ the set of non-negative integers.

\vskip 5mm

\noindent{\bf{2. Generalized Schr\"{o}dinger-Virasoro
algebras}}\vspace{3mm} \vskip 3mm

Let $\mathbb{F}$ be a field of characteristic 0, $G$ an additive
proper subgroup of $\mathbb{F}$, $\alpha\in \mathbb{F}$ satisfying
$\alpha\notin G$ while $2\alpha\in G$. We set $G_1=\alpha+G$ and
$T=G\cup G_1$. It is obvious that $T$ is an additive subgroup of
$\mathbb{F}$. In this section we want to make a natural
generalization of Schr\"{o}dinger-Virasoro algebra $\mathfrak{sv}$,
this leads us to the following definition.

\vskip 3mm

\noindent{\bf{Definition 2.1.}} The generalized
Schr\"{o}dinger-Virasoro algebra $\mathfrak{gsv}[G,\alpha]$ is
defined to be the Lie algebra with $\mathbb{F}$ basis $\{L_u, M_u,
Y_{u +\alpha}\mid u\in G\}$ subject to the following Lie brackets:
$$[L_u,L_v]=(v-u)L_{u+v},  [L_u,M_v]=v M_{u+v},
[L_u,Y_{v+\alpha}]=(v+\alpha-\frac{u}{2})Y_{u+v+\alpha},
$$
$$[M_{u},M_{v}]=0,
[Y_{u+\alpha},Y_{v+\alpha}]=(v-u)M_{u+v+2\alpha},
[M_{u},Y_{v+\alpha}]=0.$$

\vskip 3mm It is straightforward to see that
$\mathfrak{gsv}[G,\alpha]$ is $T$-graded:
$$\mathfrak{gsv}[G,\alpha]=\bigoplus\limits_{x\in
T}\mathfrak{gsv}[G,\alpha]_{x},$$ where
$$\mathfrak{gsv}[G,\alpha]_{x}=\left\{\begin{array}
{r@{\quad\quad}l} \mathbb{F}L_x\oplus\mathbb{F}M_x, & x\in G, \\
\mathbb{F}Y_x, & x\in G_1.
\end{array}\right.$$
The homogenous spaces are the root spaces according to the Cartan
subalgebra $\mathfrak{gsv}[G,\alpha]_{0}$.

One can see that if $G=\mathbb{Z}$ and $\alpha=\frac{1}{2}$, then
the generalized Schr\"{o}dinger Virasoro algebra
$\mathfrak{gsv}[G,\alpha]$ is nothing but the Schr\"{o}dinger
Virasoro algebra $\mathfrak{sv}$ defined by Henkle in [3].

\vskip 3mm Denote $L=\bigoplus\limits_{u\in G}\mathbb{F}L_u,
M=\bigoplus\limits_{u\in G}\mathbb{F}M_u, Y=\bigoplus\limits_{v\in
G}\mathbb{F}Y_{\alpha+v}, I=M\oplus Y$. Obviously, $L$ is the
centerless generalized Virasoro algebra (see [12]). $M$ and $I$ are
ideals of $L$.

\vskip 4mm

\noindent{\bf{Lemma 2.2.}} $I$ is the unique maximal ideal of
$\mathfrak{gsv}[G,\alpha]$.

\vskip 3mm

\noindent{\bf{Proof.}} It is obvious that $I$ is an ideal of
$\mathfrak{gsv}[G,\alpha]$. Moreover, $I$ is a maximal ideal of
$\mathfrak{gsv}[G,\alpha]$ since $\mathfrak{gsv}[G,\alpha]/I$ is a
simple generalized Witt algebra.

Now suppose $I_{1}$ is another maximal ideal of
$\mathfrak{gsv}[G,\alpha]$, we need to prove $I_1=I$. $\forall z\in
I_1, z\neq0$, suppose $z=l+x+y$, where $l\in L,  x\in M, y\in Y$.
Assume $l\neq 0.$ It is obvious that $0\neq adM_v(z)\in M\cap I_1$
for any $0\neq v\in G$  by using the Lie bracket of
$\mathfrak{gsv}[G,\alpha]$. Suppose
$$z_0=adM_v(z)=a_1M_{u_1}+a_2M_{u_2}+\cdots +a_nM_{u_n}.$$
It is well known that the submodules of weight module are also
weight modules, thus all homogeneous components of $z_0$ are
contained in $I_1$, hence we can find some element $0\neq u\in G$
such that $M_{u}\in I_1$. Thus $M\subseteq I_1$ due to the fact
$[L_v,M_u]=uM_{u+v}\in I_1$ for any $v\in G$. Then $$0\neq
z^{'}=l+y\in I_1.$$

If $y=0$, $z^{'}=l\in I_1$, this implies that
$I_1=\mathfrak{gsv}[G,\alpha]$, which is a contradiction. Suppose
$y\neq 0$. Since $[Y_{\alpha+v},l+y]\subseteq (M\oplus Y)\cap I_1$
for any $v\in G$ and $M\subseteq I_1$, we claim that there exists
some $0\neq y^{'}\in I_1\cap Y$. Then all homogeneous components of
$y^{'}$ are contained in $I_1$. This implies that $Y\subset I_1,
l\in I_1$, and therefore $L\subset I_1$. Therefore
$I_1=\mathfrak{gsv}[G,\alpha]$, which is a contradiction. Thus $l=0,
z=m+y$. So $I_1\subseteq M\oplus Y=I$. By the maximality of $I$ and
$I_1$, we have $I_1=I$. \hfill $\Box$

\vskip 3mm Let $G^{'}$ be another additive proper subgroup of
$\mathbb{F}$, $\alpha^{'}\in \mathbb{F}$, such that
$\alpha^{'}\notin G^{'}, 2\alpha^{'}\in G^{'}.$ Denote
$G_1^{'}=\alpha^{'}+G^{'}$, $T^{'}=G^{'}\cup G_1^{'}.$

Correspond to $G$ and $G^{'}$, there are two generalized Virasoro
algebras: $Vir[G]$ and $Vir[G^{'}]$. About $Vir[G]$ and
$Vir[G^{'}]$, one fact was pointed out in [12] that the following
Lemma can be obtained by using Theorem 4.2 in [1]. However, it can
be proved straightforward.

\vskip 3mm\noindent{\bf{Lemma 2.3.$^{[12]}$}} $Vir[G]\simeq
Vir[G^{'}]$ if and only if there exists $a\in \mathbb{F}^{*}$ such
that $aG=G^{'}$.

\vskip 3mm\noindent{\bf{Proof.}} Since $Vir[G]\simeq
Vir[G^{'}]\Leftrightarrow Vir[G]/C\simeq Vir[G^{'}]/C^{'}$, where
$C$ and $C^{'}$ are the center of $Vir[G]$ and $Vir[G^{'}]$
respectively, we view $Vir[G]$ and $Vir[G^{'}]$ as the generalized
Witt algebras. Let $\theta: Vir[G]\rightarrow Vir[G^{'}]$ be an
isomorphism of Lie algebras. Since $\mathbb{F}L_0$ and
$\mathbb{F}L_0^{'}$ are the unique Cartan subalgebras of $Vir[G]$
and $Vir[G^{'}]$ respectively, there exists $a\in \mathbb{F}^{*}$
such that $\theta(L_0)=aL_{0}^{'}$.  By using $\theta$ to
$[L_0,L_x]=xL_x$, we have
$x\theta(L_x)=[\theta(L_0),\theta(L_x)]=a[L_0^{'},\theta(L_x)],$ so
$$[L_0^{'},\theta(L_x)]=a^{-1}x\theta(L_x).$$ Thus $a^{-1}x\in
G^{'}$ since $G^{'}$ is the weight set according to the unique
Cartan algebra $\mathbb{F}L_0^{'}$. Hence $a^{-1}G\subseteq G^{'}$.
Similarly, by using $\theta^{-1}$ to $[L_0^{'},L_{x^{'}}^{'}]$, we
get $aG^{'}\subseteq G$. So, $aG^{'}=G$.

The sufficiency is obvious. This completes the proof of Lemma 2.3.

\vskip 4mm\noindent{\bf{Theorem 2.4.}} $\mathfrak{gsv}[G,\alpha]$
and $\mathfrak{gsv}[G^{'},\alpha^{'}]$ are isomorphic if and only if
there exists $a\in \mathbb{F}^{*}$ such that $G^{'}=aG, T^{'}=aT$.

\vskip 3mm

\noindent{\bf{Proof.}} Let
$\theta:\mathfrak{gsv}[G,\alpha]\rightarrow
\mathfrak{gsv}[G^{'},\alpha^{'}]$ be an isomorphism of Lie algebras,
$I$, $I^{'}$ be the maximal ideals of $\mathfrak{gsv}[G,\alpha]$ and
$\mathfrak{gsv}[G^{'},\alpha^{'}]$ respectively. By Lemma 2.2, we
have $\theta(I)=I^{'}$, thus there is an isomorphism of the
generalized Witt algebras induced by $\theta$:
$$\overline{\theta}:\mathfrak{gsv}[G,\alpha]/I\rightarrow \mathfrak{gsv}[G^{'},\alpha^{'}]/I^{'}.$$
By applying Lemma 2.3, there exists $a\in \mathbb{F}^{*}$ such that
\begin{equation} \label{eq:1}
G^{'}=aG,
\end{equation}
and
$\overline{\theta}(\overline{L}_u)=\chi(u)a^{-1}\overline{L^{'}}_{au},$
for some $\chi\in Hom(G,\mathbb{F}^{*})$. So we can assume
\begin{equation} \label{eq:1}
\theta(L_{u})=\chi(u)a^{-1}L^{'}_{au}+m_u^{'L}+y_u^{'L}.
\end{equation}

Suppose
 $\theta(Y_{\alpha+v})=\sum\limits_{i=1}^{s}c_{v_i^{'}}Y_{\alpha^{'}+v_i^{'}}^{'}+m_{v}^{'Y}.$
By using $\theta$ to both sides of the identity
$(\alpha+v)Y_{\alpha+v}=[L_0,Y_{\alpha+v}],$ we have
\begin{eqnarray*}
(\alpha+v)\theta(Y_{\alpha+v})&=&[\theta(L_0),\theta(Y_{\alpha+v})]\\
&=&[a^{-1}L_0^{'}+m_0^{'L}+y_0^{'L},\sum_{i=1}^{s}c_{v_i^{'}}Y_{\alpha^{'}+v_i^{'}}^{'}+m_{v}^{'Y}]\\
&=&[a^{-1}L_0^{'}+y_0^{'L},\sum_{i=1}^{s}c_{v_i^{'}}Y_{\alpha^{'}+v_i^{'}}^{'}+m_{v}^{'Y}]\\
&=&a^{-1}\sum\limits_{i=1}^{s}c_{v_i^{'}}[L_0^{'},Y_{\alpha^{'}+v_i^{'}}^{'}]+m{'}\\
&=&a^{-1}\sum\limits_{i=1}^{s}c_{v_i^{'}}(\alpha^{'}+v_i^{'})Y_{\alpha^{'}+v_i^{'}}^{'}+m{'}
\end{eqnarray*}
for some $m^{'}\in M^{'}$. Then we have
$$a^{-1}\sum\limits_{i=1}^{s}c_{v_i^{'}}(\alpha^{'}+v_i^{'})Y_{\alpha^{'}+v_i^{'}}^{'}+m{'}
=(\alpha+v)(\sum\limits_{i=1}^{s}c_{v_i^{'}}Y_{\alpha^{'}+v_i^{'}}^{'}+m_{v}^{'Y}).$$
By comparing the coefficients, we get
$$v_i^{'}=-\alpha^{'}+a(\alpha+v), \ \forall\ i\in \{1,\cdots, s\}.$$
Hence $s=1$ and
\begin{equation} \label{eq:1}
\theta(Y_{\alpha+v})=c_v^{'}Y_{a(\alpha+v)}^{'}+m_v^{'Y},
\end{equation}
where $c_v^{'}=c_{-\alpha^{'}+a(\alpha+v)}.$

As we know that $\theta(I)=I^{'}$, $M$ and $M^{'}$ are the centers
of $I$ and $I^{'}$ respectively, so there exists an induced
isomorphism
\begin{eqnarray*}
\bar{\bar{\theta}}: &\bar{Y}=I/M& \rightarrow  I^{'}/M^{'}=\bar{Y}^{'}:\\
         & \bar{Y}_{\alpha+v} & \mapsto
         c_v^{'}\bar{Y}_{a(\alpha+v)}^{'}.
\end{eqnarray*}
Thus we have the following isomorphisms of vector spaces:
$$Y\xrightarrow[]{\pi}\bar{Y}\xrightarrow[]{\bar{\bar{\theta}}}\bar{Y^{'}}\xrightarrow[]{\pi^{-1}}Y^{'}$$
where $\pi$ and $\pi^{'}$ are the canonical homomorphisms of vector
spaces. Thus
\begin{eqnarray*}
\pi^{'-1}\bar{\bar{\theta}}\pi: &\bigoplus\limits_{v\in G}Y_{\alpha+v}=Y& \rightarrow  Y^{'}=\bigoplus\limits_{v^{'}\in G^{'}}Y_{\alpha^{'}+v^{'}}:\\
         & Y_{\alpha+v} & \mapsto
         c_v^{'}Y_{a(\alpha+v)}^{'}
\end{eqnarray*} is an isomorphism of vector spaces. Thus
$\alpha^{'}+G^{'}=a(\alpha+G).$ This and (1) give us
$$G^{'}=aG, T^{'}=aT.$$

On the other hand, if $G^{'}=aG, T^{'}=aT$, we define a map
$$\theta: \mathfrak{gsv}[G,\alpha]\rightarrow
\mathfrak{gsv}[G^{'},\alpha^{'}]: L_u  \mapsto a^{-1}L^{'}_{au},
          M_u \mapsto a^{-1}M^{'}_{au},
          Y_{\alpha+u}\mapsto a^{-1}Y^{'}_{a(\alpha+u)}.$$
It is straightforward to check that $\theta$ is an isomorphism of
Lie algebras. This completes the proof of Theorem 2.4.\hfill $\Box$

\vskip 5mm

\noindent{\bf{Corollary 2.5.}} The map:
\begin{eqnarray*}
\theta : &\mathfrak{sv}=\mathfrak{gsv}[\mathbb{Z},\frac{1}{2}] & \rightarrow  \mathfrak{gsv}[G,\alpha]:\\
         & L_i & \mapsto (2\alpha)^{-1}L_{2\alpha i},\\
         & M_i & \mapsto (2\alpha)^{-1}M_{2\alpha i},\\
         & Y_{i+\frac{1}{2}} & \mapsto (2\alpha)^{-1}Y_{2\alpha i+\alpha}
\end{eqnarray*}extends uniquely to a Lie algebra isomorphism between
$\mathfrak{sv}=\mathfrak{gsv}[\mathbb{Z},\frac{1}{2}]$ and
$\mathfrak{gsv}[2\mathbb{Z}\alpha,\alpha]$.

\vskip 3mm \noindent{\bf{Lemma 2.6.}} Let $G, \alpha, G^{'},
\alpha^{'}$ be as in Theorem 2.4,
$\theta:\mathfrak{gsv}[G,\alpha]\rightarrow
\mathfrak{gsv}[G^{'},\alpha^{'}]$ be a Lie algebra isomorphism. Then
\vskip 3mm\begin{equation*} \label{eq:1} \left\{ \begin{aligned}
         &\theta(L_u) = \chi(u)a^{-1}L^{'}_{au}+m_u^{'L}+y_u^{'L},\\
                  &\theta(M_u)=b_uM_{au}^{'},\\
                  &\theta(Y_{\alpha+v})=c_vY^{'}_{a(\alpha+v)}+m_v^{'Y}.
                          \end{aligned} \right.
                          \end{equation*}
for some $b_u, c_v\in \mathbb{F}^{*}, m_{v}^{'Y}, m_u^{'L}\in M^{'},
y_u^{'L}\in Y^{'}$ and $\chi\in Hom(G,\mathbb{F}^{*})$, where $a\in
\mathbb{F}^{*}$ such that $aG=G^{'}, aT=T^{'}$.

\vskip 3mm \noindent{\bf{Proof.}} The first and the third identities
are given in the proof of Theorem 2.4 (i.e., identities (2) and (3)
), we only need to prove the second formula. Since
$\mathbb{F}M_0$(resp. $\mathbb{F}M_0^{'}$) is the center of
$\mathfrak{gsv}[G,\alpha]$(resp.
$\mathfrak{gsv}[G^{'},\alpha^{'}]$), we have
$$
\theta(M_0)=b_0M_0^{'}.
$$
for some $b_0\in \mathbb{F}^{*}$. Recall that $\theta(I)=I^{'}$,
$\theta(M)=M^{'}$, for $u\neq 0$, suppose
$$\theta(M_u)=\sum\limits_{i=1}^{n}b_{u_i^{'}}M_{u_i^{'}}^{'}+b_0^{'}M_{0}^{'},$$
where $u_i\neq 0, b_{u_i}\neq 0$. Then we have
\begin{eqnarray*}
u\theta(M_u)&=&\theta([L_0,M_u])\\
&=&[a^{-1}L_0^{'}+m_0^{'L}+y_0^{'L},\sum\limits_{i=1}^{n}b_{u_i^{'}}M_{u_i^{'}}^{'}+b_0^{'}M_0^{'}]\\
&=&[a^{-1}L_0^{'},\sum\limits_{i=1}^{n}b_{u_i^{'}}M_{u_i^{'}}^{'}]
=a^{-1}\sum\limits_{i=1}^{n}b_{u_i^{'}}[L_0^{'},M_{u_i^{'}}^{'}]\\
&=&a^{-1}\sum\limits_{i=1}^{n}b_{u_i^{'}}u_i^{'}M_{u_i^{'}}^{'}.
\end{eqnarray*}
So we have
$$u(\sum\limits_{i=1}^{n}b_{u_i^{'}}M_{u_i^{'}}^{'}+b_0^{'}M_{0}^{'})=a^{-1}\sum\limits_{i=1}^{n}b_{u_i^{'}}u_i^{'}M_{u_i^{'}}^{'}.$$
Thus $b_0^{'}=0,\
\sum\limits_{i=1}^{n}(ub_{u_i^{'}}-a^{-1}b_{u_i^{'}}u_i^{'})M_{u_i}^{'}=0.$
Hence $u_i^{'}=au,\ \forall i\in \{1,2,\cdots,n\},$ and $n=1$. This
gives
$$
\theta(M_u)=b_uM_{au}^{'},\ \forall u\in G.
$$
\hfill $\Box$

\vskip 5mm\noindent{\bf{3. Automorphism group of
$\mathfrak{gsv}[G,\alpha]$}}

\vskip 3mm

In this section, we first construct three kinds of automorphisms of
$\mathfrak{gsv}[G,\alpha]$ which are not inner, then determine the
automorphism group of $\mathfrak{gsv}[G,\alpha]$ completely.
Throughout this section we always assume that $\mathbb{F}$ is an
algebraic closed field with characteristic 0.

\vskip 3mm\noindent{\bf{Lemma 3.1.}} (i) For any $\chi\in
Hom(T,\mathbb{F}^{*})$ and $b\in \mathbb{F}^{*}$, the map
\begin{eqnarray*}
\sigma_{b}^{\chi} : & \mathfrak{gsv}[G,\alpha] & \rightarrow  \mathfrak{gsv}[G,\alpha]:\\
         & L_u & \mapsto \chi(u)L_{u},\\
         & M_u & \mapsto b\chi(u)M_{u},\\
         & Y_{\alpha+u} & \mapsto b^{\frac{1}{2}}\chi(\alpha+u)Y_{\alpha+u}
\end{eqnarray*}
is an automorphism of $\mathfrak{gsv}[G,\alpha].$ Furthermore, the
set $\{\sigma_{b}^{\chi}|\chi\in Hom(T,\mathbb{F}^{*}), b\in
\mathbb{F}^{*}\}$ forms a subgroup of
$Aut(\mathfrak{gsv}[G,\alpha])$ and this subgroup is isomorphic to
$(Hom(T,\mathbb{F}^{*})\times \mathbb{F}^{*})$, where
$\sigma_{b_1}^{\chi_1}\sigma_{b_2}^{\chi_2}=\sigma_{b_1b_2}^{\chi_1\chi_2}$
for $b_1, b_2\in \mathbb{F}^{*}$ and $\chi_1,\chi_2\in
Hom(T,\mathbb{F}^{*}).$

\vskip 3mm (ii) For any $a\in S(G,T):=\{a\in\mathbb{F}^{*}|aG=G,
aT=T\}$, the map
\begin{eqnarray*}
\varphi_{a}: & \mathfrak{gsv}[G,\alpha] & \rightarrow  \mathfrak{gsv}[G,\alpha]:\\
         & L_u & \mapsto a^{-1}L_{au},\\
         & M_u & \mapsto a^{-1}M_{au},\\
         & Y_{\alpha+u} & \mapsto a^{-1}Y_{a(\alpha+u)}
\end{eqnarray*}
is an automorphism of $\mathfrak{gsv}[G,\alpha],$ and
$\{\varphi_{a}|a\in S(G,T)\}$ forms a subgroup of
$Aut(\mathfrak{gsv}[G,\alpha]),$ where
$\varphi_{a}\varphi_{b}=\varphi_{ab}$ for $a, b \in \mathbb{F}^{*}.$

\vskip 3mm (iii) Denote $$\mathcal {A}:=\{\underline{a}=(a_u)_{u\in
G}|(v-u)a_{u+v}=va_v-ua_u, a_{-u}=-a_u, \forall u, v\in G\}.$$ The
map $\phi_{\underline{a}}:\mathfrak{gsv}[G,\alpha]\rightarrow
\mathfrak{gsv}[G,\alpha]$: $$L_u\mapsto L_u+a_uM_u,\ M_u\mapsto
M_u,\ Y_{\alpha+v}\mapsto Y_{\alpha+v}$$ is an automorphism.
$\Phi=\{\phi_{\underline{a}}|\underline{a}\in \mathcal {A}\}$ forms
a subgroup of $Aut(\mathfrak{gsv}[G,\alpha]),$ where
$\phi_{\underline{a}}\phi_{\underline{b}}=\phi_{\underline{a+b}}$
for $\underline{a}, \underline{b}\in \mathcal {A}.$

\vskip 3mm\noindent{\bf{Proof.}} The proof is straightforward, we
omit the details.\hfill $\Box$

\vskip 3mm \noindent{\bf{Theorem 3.2.}} Let
$Inn(\mathfrak{gsv}[G,\alpha])$ be the inner automorphism subgroup
of $Aut(\mathfrak{gsv}[G,\alpha])$. Then we have
$$Aut(\mathfrak{gsv}[G,\alpha])\simeq
(((Hom(T,\mathbb{F}^{*})\times \mathbb{F}^{*})\rtimes
S(G,T))\ltimes\Phi)\ltimes Inn(\mathfrak{gsv}[G,\alpha]).$$

\vskip 3mm \noindent{\bf{Proof.}} For any $\theta\in
Aut(\mathfrak{gsv}[G,\alpha])$, by Lemma 2.6, we can assume
\begin{eqnarray*}
         &\theta(L_u)&= \chi_1(u)a^{-1}L_{au}+m_u^{L}+y_u^{L},\\
         &\theta(M_u)&= \chi^{'}(u)a^{-1}M_{au},\\
         &\theta(Y_{\alpha+v})&= \chi^{'}(\alpha+v)a^{-1}Y_{a(\alpha+v
         )}+m_v^{Y},
\end{eqnarray*}where $\chi_1\in Hom(G,\mathbb{F}^{*}),$ $\chi^{'}: T\rightarrow
\mathbb{F}^{*}$ is a map, $a\in S(G,T)$, $m_u^{L}, m_v^{Y}\in M,
y_u^{L}\in Y.$

By applying $\theta$ to both sides of the identities:
$$[L_{-u},M_u]=uM_{0},$$
$$[L_u,Y_{\alpha+v}]=(\alpha+v-\frac{u}{2})Y_{\alpha+v+u},$$ and
$$[Y_{\alpha+u},Y_{\alpha+v}]=(v-u)M_{u+v+2\alpha},$$ we obtain
\begin{equation} \label{eq:1}
\chi^{'}(u)=b\chi_1(u),\ \forall u\in G,
\end{equation}where $b=\chi^{'}(0)\in \mathbb{F}^{*}$,
\begin{equation} \label{eq:1}
\chi_1(u)\chi^{'}(\alpha+v)=\chi^{'}(u+v+\alpha),
\end{equation}and
\begin{equation} \label{eq:1}
\chi^{'}(\alpha+u)\chi^{'}(\alpha+v)=\chi^{'}(u+v+2\alpha).
\end{equation}
for $u, v\in G$.

Write $$\chi(x)=\left\{\begin{array}
{r@{\quad\quad}l} \chi_1(x),& if\ x\in G, \\
b^{-\frac{1}{2}}\chi^{'}(x),& if\ x\in G_1=\alpha+G,
\end{array}\right.$$
and by using (4), (5) and (6), one can easily see that $\chi\in
Hom(T,\mathbb{F}^{*})$. By Lemma 3.1 (i),
\begin{eqnarray*}
\sigma_{b}^{\chi} : & \mathfrak{gsv}[G,\alpha] & \rightarrow  \mathfrak{gsv}[G,\alpha]:\\
         & L_u & \mapsto \chi(u)L_{u},\\
         & M_u & \mapsto b\chi(u)M_{u},\\
         & Y_{\alpha+v} & \mapsto b^{\frac{1}{2}}\chi(\alpha+v)Y_{\alpha+v}
\end{eqnarray*}
is an automorphism of $\mathfrak{gsv}[G,\alpha]$.

For $(\sigma_{b}^{\chi})^{-1}\theta\in
Aut(\mathfrak{gsv}[G,\alpha])$, it is obvious that
\begin{eqnarray*}
(\sigma_{b}^{\chi})^{-1}\theta: & \mathfrak{gsv}[G,\alpha] & \rightarrow  \mathfrak{gsv}[G,\alpha]:\\
         & L_u & \mapsto a^{-1}L_{au}+m_{1u}^{L}+y_{1u}^{L},\\
         & M_u & \mapsto a^{-1}M_{au},\\
         & Y_{\alpha+v} & \mapsto a^{-1}Y_{a(\alpha+v)}+m_{1v}^{Y},
\end{eqnarray*}
where $m_{1u}^{L},m_{1v}^{Y}\in M$, $y_{1u}^{L}\in Y$.

Set
\begin{eqnarray*}
\varphi_{a}: & \mathfrak{gsv}[G,\alpha] & \rightarrow  \mathfrak{gsv}[G,\alpha]:\\
         & L_u & \mapsto a^{-1}L_{au},\\
         & M_u & \mapsto a^{-1}M_{au},\\
         & Y_{\alpha+v} & \mapsto a^{-1}Y_{a(\alpha+v)}.
\end{eqnarray*}
By Lemma 3.1 (ii), $\varphi_{a}\in Aut(\mathfrak{gsv}[G,\alpha])$.
Set $\tau=(\varphi_{a})^{-1}(\sigma_{b}^{\chi})^{-1}\theta$, then
$\tau\in Aut(\mathfrak{gsv}[G,\alpha]).$ More precisely
\begin{equation} \label{eq:1}
\tau(L_u)=L_{u}+m_{2u}^{L}+y_{2u}^{L},\\
\tau(M_u)=M_{u},\\
\tau(Y_{\alpha+v})=Y_{\alpha+v}+m_{2v}^{Y},
\end{equation}
where $m_{2u}^{L},m_{2v}^{Y}\in M$, $y_{2u}^{L}\in Y$.

\vskip3mm \noindent{\bf{Claim.}}
$\tau=(\varphi_{a})^{-1}(\sigma_{b}^{\chi})^{-1}\theta\in
Inn(\mathfrak{gsv}[G,\alpha])\cdot\Phi$.

\vskip 3mm In fact, assume
$$\tau(L_0)=L_0+\sum\limits_{i=1}^{p}a_{u_i}M_{u_i}+\sum\limits_{j=1}^{q}b_{v_j}Y_{\alpha+v_{j}}+b_0Y_{\alpha},\
\tau(Y_{\alpha})=Y_{\alpha}+\sum\limits_{k=1}^{r}c_{w_k}M_{w_k}.$$
Applying $\tau$ to $[L_0,Y_{\alpha}]=\alpha Y_{\alpha}$, we have
$$\sum\limits_{k=1}^{r}c_{w_k}w_kM_{w_k}-\sum\limits_{j=1}^{q}b_{v_j}v_jM_{2\alpha+v_j}=\alpha\sum\limits_{k=1}^{r}
c_{w_k}M_{w_k}.$$ Noting that $v_j\neq 0,\ \alpha\neq w_k$, we have
$$q=r,\ v_k=w_k-2\alpha,\
b_{v_k}=\frac{c_{w_k(w_k-\alpha)}}{w_k-2\alpha}.$$ These gives
us
\begin{equation} \label{eq:1}
\tau(L_0)=L_0+\sum\limits_{i=1}^{p}a_{u_i}M_{u_i}+\sum\limits_{k=1}^{r}
\frac{c_{w_k}(w_k-\alpha)}{w_k-2\alpha}Y_{w_k-\alpha}+b_0Y_{\alpha}.
\end{equation}
Now we construct an inner automorphism $tau^{'}$ for
$gsv[G,\alpha]$, which is equal to $tau$ when acting on $L_0$ and
$Y_{\alpha}$. Indeed, we set
\begin{eqnarray*}
\tau^{'}&=&expad(\sum\limits_{1\leq j\neq k\leq
r}\frac{c_{w_k}c_{w_j}(\alpha-w_j)(w_j-w_k)}
{2(w_k-2\alpha)(w_j-2\alpha)(w_j+w_k-2\alpha)}M_{w_j+w_k-2\alpha}\\
& &
+\sum\limits_{k=1}^{r}\frac{b_0c_{w_k}(2\alpha-w_k)}{2\alpha w_k}M_{w_k})\\
& & expad(-\sum\limits_{i=1}^{p}\frac{a_{u_i}}{u_i}M_{u_i})expad
(\sum\limits_{k=1}^{r}\frac{c_{w_k}}{2\alpha-w_k}Y_{w_k-\alpha}-\frac{b_0}{\alpha}Y_{\alpha}),
\end{eqnarray*}
Then one can check that the inner automorphism $\tau^{'}$ satisfying
$$\tau^{'}(L_0)=\tau(L_0),  \tau^{'}(Y_{\alpha})=\tau(Y_{\alpha}).$$

For any $u\in G$, we apply $\tau^{'}$ to $[L_0,L_u]=uL_{u}$. For the
right side, we have
$$
u\tau^{'}(L_u)=u\tau(L_u)+u(\tau^{'}(L_u)-\tau(L_u)).
$$
For the left side, we have
\begin{equation*} \label{eq:3}
\begin{split}
[\tau^{'}(L_0),\tau^{'}(L_u)]&= [\tau(L_0),\tau^{'}(L_u)]=
[\tau(L_0),\tau(L_u)+(\tau^{'}(L_u)-\tau(L_u))] \\
 &= u\tau(L_u)+[\tau(L_0),\tau^{'}(L_u)-\tau(L_u)].
 \end{split}
 \end{equation*}
Thus
\begin{equation} \label{eq:1}
[\tau(L_0),\tau^{'}(L_u)-\tau(L_u)]=u(\tau^{'}(L_u)-\tau(L_u)).
\end{equation}
Now we prove the following identity.
\begin{equation} \label{eq:1}
\tau^{'}(L_u)=\tau(L_u)+e_uM_u.
\end{equation}
for some $e_u\in\mathbb{F}$.

Indeed, By using (7) and the definition of $\tau^{'}$, one can see
that
\begin{equation} \label{eq:1} \tau^{'}(M_u)=\tau(M_u)=M_u.
\end{equation} So
$[\tau^{'}(L_u),M_v]=vM_{u+v}.$ Thus $\tau^{'}(L_u)=L_u+m_1+y_1$ for
some $m_1\in M, y_1\in Y.$ By using (7) again, we have
$\tau^{'}(L_u)-\tau(L_u)\in M\oplus Y.$ Assuming
$$\tau^{'}(L_u)-\tau(L_u)=\sum_{k=1}^{r}e_{v_k}M_{v_k}+\sum_{l=1}^{s}d_{w_l}Y_{\alpha+w_l}.$$
By using this along with identities (8) and (9), we have
\begin{eqnarray*}
&&[L_0+\sum\limits_{i=1}^{p}a_{u_i}M_{u_i}+\sum\limits_{j=1}^{q}b_{v_j}Y_{\alpha+v_{j}}+b_0Y_{\alpha},
\sum_{k=1}^{r}e_{v_k}M_{v_k}+\sum_{l=1}^{s}d_{w_l}Y_{\alpha+w_l}]\\
&=&\sum_{k=1}^{r}e_{v_k}v_kM_{v_k}+\sum_{l=1}^{s}d_{w_l}(\alpha+w_l)Y_{\alpha+w_l}\\
&&+\sum_{j=1}^{q}\sum_{l=1}^{s}b_{v_j}d_{w_l}(w_l-v_j)M_{2\alpha+w_l+v_j}
+\sum_{l=1}^{s}b_0d_{w_l}w_lM_{2\alpha+w_l}\\
&=&u(\sum_{k=1}^{r}e_{v_k}M_{v_k}+\sum_{l=1}^{s}d_{w_l}Y_{\alpha+w_l}).
\end{eqnarray*}
By comparing the coefficients of $Y_{\alpha+w_l}$ we have
$\sum_{l=1}^{s}ud_{w_l}=\sum_{l=1}^{s}d_{w_l}(\alpha+w_l)$. This
means $d_{w_l}=0$ for any $l\in \{1,\cdots,s\}$ since $u\neq
\alpha+w_l.$ Thus the we get
$$\sum_{k=1}^{r}e_{v_k}v_kM_{v_k}=u(\sum_{k=1}^{r}e_{v_k}M_{v_k}).$$
Furthermore, $r=1, \ u=v_1$. Thus $\tau^{'}(L_u)-\tau(L_u)=e_uM_u.$
This proves (10).

For any $v\in G, v\neq 2\alpha$, by (10), we have
$$\tau^{'}(Y_{\alpha+v})=(\alpha-\frac{v}{2})^{-1}\tau^{'}[L_v,Y_{\alpha}]=
(\alpha-\frac{v}{2})^{-1}[\tau(L_v)+e_vM_v,\tau(Y_{\alpha})]=\tau(Y_{\alpha+v}).$$
For $v=2\alpha$, we have
$$\tau^{'}(Y_{3\alpha})=(-3\alpha)^{-1}[\tau^{'}(L_{4\alpha}),\tau^{'}(Y_{-\alpha})]=
(-3\alpha)^{-1}[\tau(L_{4\alpha})+b_{4\alpha}M_{4\alpha},\tau(Y_{-\alpha})]=\tau(Y_{3\alpha}).$$
In all cases we get
\begin{equation} \label{eq:1}
\tau^{'}(Y_{\alpha+v})=\tau(Y_{\alpha+v}),\ \forall \ v\in G.
\end{equation}

 Define$$\phi:\mathfrak{gsv}[G,\alpha]\rightarrow
\mathfrak{gsv}[G,\alpha]: L_u\mapsto L_u+e_uM_u,\ M_u\mapsto M_u,\
Y_{\alpha+v}\mapsto Y_{\alpha+v}.$$ By (10), (11) and (12), we have
$\tau^{'}=\tau\phi.$ Thus $\tau=\tau^{'}\phi^{-1}\in
Inn(\mathfrak{gsv}[G,\alpha])\cdot\Phi.$ Which completes the proof
of the Claim.

\vskip 3mm By the claim, we have
$$\theta=\sigma_b^{\chi}\varphi_a\tau^{'}\phi^{-1}\in
(Hom(T,\mathbb{F}^{*})\times \mathbb{F}^{*})\cdot S(G,T)\cdot
Inn(\mathfrak{gsv}[G,\alpha])\cdot \Phi.$$ Since
$Inn(\mathfrak{gsv}[G,\alpha])$ is a normal subgroup of
$Aut(\mathfrak{gsv}[G,\alpha])$, thus
$Inn(\mathfrak{gsv}[G,\alpha])\cdot \Phi=\Phi\cdot
Inn(\mathfrak{gsv}[G,\alpha])$. One can check straightforward that
the following two facts hold:
$$(Hom(T,\mathbb{F}^{*})\times \mathbb{F}^{*})\triangleleft (Hom(T,\mathbb{F}^{*})\times \mathbb{F}^{*})\cdot S(G,T),$$
and $$\Phi\triangleleft(Hom(T,\mathbb{F}^{*})\times
\mathbb{F}^{*})\cdot S(G,T)\cdot \Phi.$$ Thus
$$Aut(\mathfrak{gsv}[G,\alpha])\simeq
(((Hom(T,\mathbb{F}^{*})\times \mathbb{F}^{*})\rtimes
S(G,T))\ltimes\Phi)\ltimes Inn(\mathfrak{gsv}[G,\alpha]).$$ This
completes the proof of the Theorem 3.2.\hfill $\Box$

 \vskip 5mm

\noindent{\bf{4. Verma modules of $\mathfrak{gsv}[G,\alpha]$}}

\vskip 3mm

In this section we construct and investigate the structure of Verma
modules over the generalized Schr\"{o}dinger-Virasoro algebra
$\mathfrak{gsv}[G,\alpha]$.

Note that $T=G\cup G_1$ is a subgroup of $\mathbb{F}$, we fix a
total order $"\succeq"$ on $T$ which is compatible with the
addition, i.e., $x\succeq y$ implies $x+z\succeq y+z$ for any $z\in
T$ (see [4],[9]). We write $x\succ y$ if $x\succeq y$ and $x\neq y$.
Let
$$T_{+}:=\{x\in T|x\succ 0\},\ T_{-}:=\{x\in T|x\prec 0\}.$$ Then
$T=T_{+}\cup\{0\}\cup T_{-}$ and $\mathfrak{gsv}[G,\alpha]$ has a
triangular decomposition:
$$\mathfrak{gsv}[G,\alpha]=\mathfrak{gsv}[G,\alpha]_{-}\oplus\mathfrak{gsv}[G,\alpha]_{0}\oplus\mathfrak{gsv}[G,\alpha]_{+},$$
where $$\mathfrak{gsv}[G,\alpha]_{-}=\bigoplus_{u\prec
0}\mathbb{F}L_u\oplus \bigoplus_{u\prec 0}\mathbb{F}M_u\oplus
\bigoplus_{\alpha+v\prec 0}\mathbb{F}Y_{\alpha+v},$$
$$\mathfrak{gsv}[G,\alpha]_{+}=\bigoplus_{u\succ
0}\mathbb{F}L_u\oplus \bigoplus_{u\succ 0}\mathbb{F}M_u\oplus
\bigoplus_{\alpha+v\succ 0}\mathbb{F}Y_{\alpha+v}$$ and
$\mathfrak{gsv}[G,\alpha]_{0}=\mathbb{F}L_0\oplus \mathbb{F}M_0.$
The universal enveloping algebra of $\mathfrak{gsv}[G,\alpha]$ is
given by
$$U(\mathfrak{gsv}[G,\alpha])=U(\mathfrak{gsv}[G,\alpha])_{-}U(\mathfrak{gsv}[G,\alpha])_{0}U(\mathfrak{gsv}[G,\alpha])_{+}.$$
The elements $L_{i_1}\cdots L_{i_r}M_{j_1}\cdots
M_{j_s}Y_{\alpha+k_1}\cdots Y_{\alpha+k_t},$ where $r,s,t\in
\mathbb{N}, i_1\succeq \cdots \succeq i_{r}, j_1\succeq \cdots
\succeq j_{s}, k_1\succeq \cdots \succeq k_{t},$ along with $1$,
form a basis of $U(\mathfrak{gsv}[G,\alpha])$.

Let $c, h\in \mathbb{F}$, $V_h$ be a 1-dimensional vector space over
$\mathbb{F}$ spanned by $v_h$, i.e., $V_h=\mathbb{F}v_h$. View $V_h$
as a $\mathfrak{gsv}[G,\alpha]_0$-module such that $L_0.v_h=hv_v,$
$M_{0}.v_h=cv_h$. Then $V_h$ is a
$\mathfrak{B}=\mathfrak{gsv}[G,\alpha]_{+}\oplus\mathfrak{gsv}[G,\alpha]_{0}$-module
by setting $\mathfrak{gsv}[G,\alpha]_{+}.V_h=0$.

\vskip 3mm

\noindent{\bf{Definition 4.1.}} The induced module
$V(c,h)=Ind_{\mathfrak{B}}^{\mathfrak{gsv}[G,\alpha]}V_h=U(\mathfrak{gsv}[G,\alpha])\otimes_{U(\mathfrak{B})}V_h$
is called the Verma module of $\mathfrak{gsv}[G,\alpha]$ with
highest weight $(c,h)$.

\vskip 3mm Let $U:=U(\mathfrak{gsv}[G,\alpha]).$ For any $c, h\in
\mathbb{F}$, let $I(c,h)$ be the left ideal of $U$ generated by the
elements
$$\{L_u, M_u, Y_{\alpha+v}|u\in G_{+}, \alpha+v\in G_{1+}\}\cup \{L_0-h, M_0-c\},$$
where $G_{+}=G\cap T_{+}, G_{1+}=G_1\cap T_{+}.$ Then the Verma
module with highest weight $(c,h)$ for $\mathfrak{gsv}[G,\alpha]$
also can be defined as $V(c,h):=U/I(c,h).$

By definition, we can easily get a basis of $V(c,h)$ consisting of
all vectors of the form:
$$v_h, L_{-i_1}\cdots L_{-i_r}M_{-j_1}\cdots
M_{-j_s}Y_{\alpha-k_1}\cdots Y_{\alpha-k_t}v_h,$$ where $$ 0\prec
i_1\preceq\cdots \preceq i_{r}, 0\prec j_1\preceq \cdots \preceq
j_{s},\alpha\prec k_1\preceq \cdots \preceq k_{t};r,s,t\in
\mathbb{N}.$$

\vskip 3mm\noindent{\bf{Remark.}} One can see that $M_0$ acts as a
scalar $c$ on $V(c,h)$ since $\mathbb{F}M_0$ is the center of
$\mathfrak{gsv}[G,\alpha]$. Next, we call a vector $v\in V(c,h)$ a
weight vector with weight $\mu$ means $v$ satisfying $L_0v=\mu v$.

\vskip 3mm\noindent{\bf{Lemma 4.2.}} $V(c,h)$ is a weight module of
$\mathfrak{gsv}[G,\alpha]$, and $V(c,h)=\oplus_{\mu\in h-T_+}
V_\mu$, where $V_\mu=span\{L_{-i_1}\cdots L_{-i_r}M_{-j_1}\cdots
M_{-j_s}Y_{\alpha-k_1}\cdots Y_{\alpha-k_t}v_h\mid
\sum_{p=1}^{r}i_p-\sum_{p=1}^{s}j_p-\sum_{p=1}^{t}(\alpha-k_p)=h-\mu\}$
is the weight vector space with weight $\mu$.

\vskip 3mm\noindent{\bf{Proof.}}  It suffices to show that $L_0$
acts diagonally on the basis elements of $V(c,h)$. By the definition
of $v_h$, $L_0v_h=hv_h$. Suppose $u\in V(c,h)$ such that $L_0u=au$.
Then
$$L_0(L_{-i}u)=(a-i)L_{-i}u, L_0(M_{-j}u)=(a-j)M_{-j}u,
L_0(Y_{\alpha-k}u)=(a+\alpha-k)Y_{\alpha-k}u.$$ Thus Lemma 4.2
holds.\hfill $\Box$

\vskip 3mm We know from [4] that for the fixed total order
$"\succeq"$ of $T$, either $"\succeq"$ is dense, i.e., $\forall x\in
T_{+},$ the cardinality of $\{y\in T|0\prec y\prec x\}$ is infinite,
or $"\succeq"$ is discrete, i.e., there exists $a\in T$ such that
the set $\{y\in T|0\prec y\prec a\}$ is empty.

For the generalized Virasoro algebra $Vir[G]$ studied in [4], the
irreducibility of Verma module over $Vir[G]$ is depends on whether
the total order of $G$ is dense or discrete (see Theorem 3.1 in
[4]). With respect to Verma modules over generalized Witt algebras
studied in [8], the irreducibility depends on the action of $L_0$ on
the highest weight vector (see Theorem 3 in [9]). It is very
interesting that the irreducibility of Verma modules over
$\mathfrak{gsv}[G,\alpha]$ depends on neither the action of $L_0$
nor whether the total order is dense or discrete, we point out that
the irreducibility just depends on the action of the element $M_0$.

\vskip 5mm For $x\in V(c,h)$, we set
$$
x=\sum\limits_{\begin{subarray} \ \ \ \ \ i_1\preceq \cdots \preceq
i_{r},j_1\preceq \cdots \preceq j_{s}, k_1\preceq \cdots \preceq
k_{t}\\ i_1,\cdots, i_r,j_1,\cdots, j_s\in G_{+},k_1-\alpha,\cdots,
k_t-\alpha\in G_{1+}
\end{subarray}}
a_{\underline{i},\underline{j},\underline{k}}L_{-i_1}\cdots
L_{-i_r}M_{-j_1}\cdots M_{-j_s}Y_{\alpha-k_1}\cdots
Y_{\alpha-k_t}v_h.
$$
where $a_{\underline{i},\underline{j},\underline{k}}\in
\mathbb{F},\underline{i}=(i_1,\cdots,i_r),
\underline{j}=(j_1,\cdots,j_s),
\underline{k}=(\alpha-k_1,\cdots,\alpha-k_t)$, and only finitely
many $a_{\underline{i},\underline{j},\underline{k}}\neq 0$. We
define $$A_x:=\{\underline{i}=(i_1,\cdots,i_r)|
a_{\underline{i},\underline{j},\underline{k}}\neq 0 \ for\ some\
\underline{j},\underline{k}\},l=max\{r|\underline{i}=(i_1,\cdots,i_r)\in
A_x\},$$ where $l=0$ if $A_x=\emptyset$. We also define $l$ to be
the length of the element $x$, and denote it by $len(x)$, i.e.,
$l=len(x)$.

For $r\in \mathbb{N}$, we set
$$V_{r}:=span_{\mathbb{F}}\{x|len(x)\leq r\}. $$ In what follows, we assume
$V_r=0$ if $r\leq -1$. One can check the following two lemmas by
straightforward and easy computations.

\vskip 3mm\noindent{\bf{Lemma 4.3.}} (i)
\begin{eqnarray*}
& &M_jL_{-i_1}L_{-i_2}\cdots L_{-i_r}M_{-j_1}\cdots
M_{-j_s}Y_{\alpha-k_1}\cdots Y_{\alpha-k_t}v_h\\
&\equiv&-j(\sum_{1\leq p\leq r}L_{-i_1}\cdots \hat{L}_{-i_p}\cdots
L_{-i_r}M_{j-i_p})M_{-j_1}\cdots M_{-j_s}Y_{\alpha-k_1}\cdots
Y_{\alpha-k_t}v_h(modV_{r-2}),
\end{eqnarray*}
for any $r\in \mathbb{N}; j\in G_{+}; 0\prec i_1\preceq\cdots\preceq
i_r; 0\prec j_1\preceq\cdots\preceq j_s; 0\prec
k_1-\alpha\preceq\cdots\preceq k_t-\alpha,$ where $\hat{  }$ means
the corresponding element is deleted.

\vskip 3mm(ii)
\begin{eqnarray*}
& &L_jM_{-i_1}M_{-i_2}\cdots M_{-i_r}v_h\\
&=& (\sum_{1\leq p\leq r}(-i_p)M_{-i_1}\cdots \hat{M}_{-i_p}\cdots
M_{-i_r} M_{j-i_p})v_h,
\end{eqnarray*}
for any $r\in \mathbb{N}; j,i_1,\cdots,i_r\in G_{+}.$ In
particularly, $$L_jM_{-i_1}M_{-i_2}\cdots M_{-i_r}v_h=0, \forall
j\succ max\{i_1,\cdots, i_r\}.$$

\vskip 3mm\noindent{\bf{Lemma 4.4.}}
(i)$Y_{-\alpha+j}Y_{\alpha-k_1}Y_{\alpha-k_2}\cdots
Y_{\alpha-k_t}v_h=0,\ \forall j\succ k_t,$  where $\alpha\prec
k_1\preceq \cdots\preceq k_t$.

\vskip 3mm(ii)If $M_0.v_h=0$, then
$Y_{-\alpha+j}Y_{\alpha-k_1}Y_{\alpha-k_2}\cdots
Y_{\alpha-k_t}v_h=0, \ \forall j\succeq k_t,$ where $\alpha\prec
k_1\preceq \cdots\preceq k_t$.

\vskip 3mm\noindent{\bf{Corollary 4.5.}} $M_jV_r\subseteq V_{r-1}$,
 for any $j\in G_{+}$.

\vskip 3mm \noindent{\bf{Proof.}} It follows Lemma 4.3 (i)
immediately.\hfill $\Box$

\vskip 5mm \noindent{\bf{Theorem 4.6.}} (i) The Verma module
$V(c,h)$ is an irreducible $\mathfrak{gsv}[G,\alpha]$ module if
$c\neq 0$. \vskip 3mm

(ii) If $c=0$, then the Verma module $V(0,h)$ contains a unique
maximal proper submodule $N(0,h)$, where $N(0,h)$ is generated by
$\{L_{-u}v_h, M_{-u}v_h, Y_{\alpha-v}v_h|u\in G_+, v-\alpha\in
G_{1+}\}$ if $h=0$, by $\{M_{-u}v_h, Y_{\alpha-v}v_h|u\in G_+,
v-\alpha\in G_{1+}\}$ if $h\neq 0$.

\vskip 3mm

\noindent{\bf{Proof.}} (i) Suppose $c\neq 0$. Let $u_0\neq 0$ be any
given weight vector in $V(c,h)$. By Lemma 4.2 and the fact that a
submodule of a weight module is a weight module, we need only to
prove that $v_h\in U(\mathfrak{gsv}[G,\alpha])u_0$.

\vskip 3mm

\noindent{\bf{Claim I.}} There exists a weight vector $u\in
U(\mathfrak{gsv}[G,\alpha])u_0$ such that
$$
u=\sum\limits_{\begin{subarray}\  \ \ \ \ j_1\preceq \cdots \preceq
j_{s}; \  k_1\preceq \cdots \preceq k_{t}\\  j_1,\cdots, j_s\in
G_{+}; k_1-\alpha,\cdots, k_t-\alpha\in G_{1+}
\end{subarray}}a_{\underline{j},\underline{k}}
M_{-j_1}\cdots M_{-j_s}Y_{-k_1+\alpha}\cdots Y_{-k_t+\alpha}v_h,
$$
where $a_{\underline{j},\underline{k}}\in \mathbb{F}$ and only
finitely many $a_{\underline{j},\underline{k}}\neq 0$,
$\underline{j}=(j_1,\cdots,j_s), \underline{k}=(k_1,\cdots,k_t).$

\vskip 3mm

In fact, suppose
$$
u_0=\sum\limits_{\begin{subarray}\  0\prec  i_1\preceq
\cdots \preceq i_{r},\\ 0\prec j_1\preceq \cdots \preceq j_{s},\\
\alpha\prec k_1\preceq \cdots \preceq k_{t}.
\end{subarray}}
a_{\underline{i},\underline{j},\underline{k}}L_{-i_1}\cdots
L_{-i_r}M_{-j_1}\cdots M_{-j_s}Y_{\alpha-k_1}\cdots
Y_{\alpha-k_t}v_h
$$
$$
\equiv\sum\limits_{\begin{subarray}\  0\prec i_1\preceq
\cdots \preceq i_{l},\\ 0\prec j_1\preceq \cdots \preceq j_{s},\\
\alpha\prec k_1\preceq \cdots \preceq k_{t}.
\end{subarray}}
a_{\underline{i},\underline{j},\underline{k}}L_{-i_1}\cdots
L_{-i_l}M_{-j_1}\cdots M_{-j_s}Y_{\alpha-k_1}\cdots
Y_{\alpha-k_t}v_h(modV_{l-1}),
$$
where $l=len(u_0), \underline{i}=(i_1,\cdots,i_r), \underline{j},
\underline{k}$ as above,
$a_{\underline{i},\underline{j},\underline{k}}\in \mathbb{F}$.

If $l=0$, there is nothing to prove. Now suppose $len(u_0)=l\geq 1$,
we denote
$$i_{l}^{(0)}:=max\{i_l|(i_1,\cdots,i_l)\in A_{u_0}\},$$ where $A_{u_0}:=
\{\underline{i}|\underline{i}=(i_1,\cdots,i_l),
a_{\underline{i},\underline{j},\underline{k}}\neq 0 \ for\ some\
\underline{j},\underline{k}\}$, then by using Lemma 4.3 (i) and
Corollary 4.4 we can deduce that
$$
u_1 = M_{i_{l}^{(0)}}u_0=\sum
a_{\underline{i},\underline{j},\underline{k}}M_{i_{l}^{(0)}}L_{-i_1}\cdots
L_{-i_r}M_{-j_1}\cdots M_{-j_s}Y_{\alpha-k_1}\cdots
Y_{\alpha-k_t}v_h
$$
$$ \equiv \sum
a_{\underline{i}^{(1)},\underline{j}^{(1)},\underline{k}^{(1)}}L_{-i_1^{(1)}}\cdots
L_{-i_{l-1}^{(1)}}M_{-j_1^{(1)}}\cdots
M_{-j_{s^{'}}^{(1)}}Y_{\alpha-k_1^{(1)}}\cdots
Y_{\alpha-k_{t^{'}}^{(1)}}v_h(modV_{l-2}).
$$
it is clear that $u_1\neq 0$ and $len(u_1)=l-1$.

Repeating the precess and define $u_s$ recursively for
$s=2,\cdots,l$, one obtains the claim.

\vskip 3mm

\noindent{\bf{Claim II.}} There exists an weight vector $w\in
U(\mathfrak{gsv}[G,\alpha])u_0$ such that $w$ takes the following
form
$$w=\sum\limits_{\begin{subarray} \ \ s\geq 0,
j_1\preceq \cdots \preceq j_{s}; j_1,\cdots, j_s\in G_{+}
\end{subarray}}
a_{\underline{j}}M_{-j_1}\cdots M_{-j_s}v_h.$$

 In fact, by Claim I, we know that there is a weight vector $u\in
U(\mathfrak{gsv}[G,\alpha])u_0$ such that
$$
u=\sum\limits_{\begin{subarray}\ \ \ \ \ \ \  s,t\geq 0, j_1\preceq
\cdots \preceq j_{s}\\ \ \ \  k_1\preceq \cdots \preceq k_{t};
j_1,\cdots, j_s\in G_{+}\\ \ \ \ \ \ k_1-\alpha,\cdots,
k_t-\alpha\in G_{1+}
\end{subarray}}a_{\underline{j},\underline{k}}
M_{-j_1}\cdots M_{-j_s}Y_{-k_1+\alpha}\cdots Y_{-k_t+\alpha}v_h.
$$
  Set $$B:=\{\underline{k}=(k_1,k_2,\cdots,k_t)|a_{\underline{j},\underline{k}}\neq 0\ for\ some\ \underline{j}\},
  k^{(0)}=max\{k_t|\underline{k}\in B\}.$$ Then by Lemma 4.4, we have

$$w_1=Y_{-(\alpha-k^{(0)})}u=\sum a_{\underline{j},\underline{k}^{(1)}}M_{-j_1}\cdots
M_{-j_s}Y_{\alpha-k_1^{(1)}}\cdots
Y_{\alpha-k_{t-1}^{(1)}}(M_0v_h).$$

Noting that $w_1\neq 0$ and $w_1\in U(\mathfrak{gsv}[G,\alpha])u_0$
is a weight vector. One repeats the precess to get Claim II.

\vskip 3mm From Claim II we know that there exists a weight vector
$w\in U(\mathfrak{gsv}[G,\alpha])u_0$ such that it has the following
form
$$w=\sum\limits_{\begin{subarray} \ \ s\geq 0, 0\prec
j_1\preceq \cdots \preceq j_{s}
\end{subarray}}
a_{\underline{j}}M_{-j_1}\cdots M_{-j_s}v_h.$$  We define
$$length(w)=max\{s|\underline{j}=(j_1,\cdots,j_s),a_{\underline{j}}\neq
0\}.$$ If $length(w)=0$, then $v_h\in
U(\mathfrak{gsv}[G,\alpha])u_0$ and (i) holds. Now suppose
$length(w)>0$. Denote
$j^{(0)}=max\{j_s|\underline{j}=(j_1,\cdots,j_s),a_{\underline{j}}\neq
0\}$. By applying $L_{j^{(0)}}$ to $w$ and using Lemma 4.3 (ii), we
have
$$0\neq w_1=L_{j^{(0)}}w=\sum\limits_{\begin{subarray} \ \
j_1^{(1)}\preceq \cdots \preceq j_{s}^{(1)}; j_1^{(1)},\cdots,
j_s^{(1)}\in G_{+}
\end{subarray}}
a_{\underline{j}}^{(1)}M_{-j_1^{(1)}}\cdots M_{-j_s^{(1)}}M_0v_h.$$
It is clear that $$length(w_1)<length(w).$$ Repeating the precess,
we obtain$$0\neq w_s=aM_0v_h=acv_h\in
U(\mathfrak{gsv}[G,\alpha])u_0$$ for some $0\neq a\in \mathbb{F}.$
So $v_h\in U(\mathfrak{gsv}[G,\alpha])u_0$ and $V(c,h)$ is
irreducible.

\vskip 3mm

(ii) If $c=0, h=0$, by the definition of $N(0,0)$, one knows that
all the basis elements of $V(0,0)$ except $v_h$ are clearly in
$N(0,0)$. It suffices to show that $v_h\notin N(0,0)$. For any
weight vector $v\in N(0,0)$, suppose the weight of $v$ is $\mu$, and
for any basis element $L_{i_1}\cdots L_{i_r}M_{j_1}\cdots
M_{j_s}Y_{\alpha+k_1}\cdots Y_{\alpha+k_t}$ of
$U(\mathfrak{gsv}[G,\alpha])$ such that
$\sum_{p=1}^{r}i_p+\sum_{p=1}^{s}j_p+\sum_{p=1}^{t}(\alpha+k_p)=-\mu$,
we have $$L_{i_1}\cdots L_{i_r}M_{j_1}\cdots
M_{j_s}Y_{\alpha+k_1}\cdots Y_{\alpha+k_t}v=aL_0v_h+bM_0v_h=0,$$ for
some $a, b\in \mathbb{F}$. This implies that $v_h\notin N(0,0)$.

\vskip 3mm If $c=0, h\neq 0$, similarly as above, we can see that
$U(L_-)\notin N(0,h)$, where $L_-=\oplus_{u\prec 0}\mathbb{F}L_u$.
This means that $N(0,h)$ is a proper submodule of $V(0,h)$. Suppose
$V$ is any submodule of $V(0,h)$ such that $V\varsupsetneq N(0,h)$,
then there exist $i_1,\cdots, i_r\in G_{+}, r\in \mathbb{N}$ such
that $L_{-i_1}\cdots L_{-i_r}v_h\in V.$ If $r=0$, then $v_h\in V$
and $V=V(0,h)$. Suppose $r\geq 1$. We denote $i=i_1+i_2+\cdots
+i_r$, then
$$L_iL_{-i_1}\cdots L_{-i_r}v_h=(-1)^r(i+i_1)(i-i_1+i_2)\cdots
(i-i_1-i_2-\cdots -i_{r-1}+i_r)hv_h\in V.$$  Since
$(-1)^r(i+i_1)(i-i_1+i_2)\cdots (i-i_1-i_2-\cdots -i_{r-1}+i_r)h\neq
0$ we have $v_h\in V$ and $V=V(0,h)$. So $N(0,h)$ is the unique
maximal proper submodule of $V(0,h)$.\hfill $\Box$

\vskip 3mm\noindent{\bf{Remark.}}  $V(0,0)/N(0,0)\simeq \mathbb{F}$
is a trivial module of $\mathfrak{gsv}[G,\alpha]$.
$V(0,h)/N(0,h)\simeq U(L_-)$ as vector space.

\vskip 3mm\noindent{\bf{Acknowledgements.}} We would like to thank
Prof. Kaiming Zhao for providing the proof of Lemma 2.3.

\end{document}